\documentclass{amsart}
\usepackage[latin1]{inputenc}
\usepackage{amsmath,amssymb,amscd,latexsym}
\usepackage[english]{babel}
\usepackage{color,epsfig}

\usepackage{pst-tree}
\newcommand{\arbre}[2]{
\pstree[levelsep=0.7cm, xbbl=0.3cm, xbbr=0.3cm, xbbh=0.3cm,
xbbd=0.3cm]{\noeud{#1}}{#2}}
\newcommand{\noeud}[1]{%
   \TR{\rule[-0.1cm]{0mm}{0.1cm}#1}%
   }
\newtheorem{theorem}{Theorem}[section]
\newtheorem{lemma}[theorem]{Lemma}
\newtheorem{defi}[theorem]{Definition}
\newtheorem{exam}[theorem]{Example}
\newtheorem{cons}[theorem]{Construction}

\newcommand{\Pn}{\mathbb P_n}

\newcommand{\Pnt}{\tilde{\mathbb P}_n}
\newcommand{\Pnb}{\bar{\mathbb P}_n}
\renewcommand{\P}{\mathbb P}
\newcommand{\Q}{\mathbb Q}
\newcommand{\R}{\mathbb R}

\newcommand{\E}{\mathbb E}
\newcommand{\N}{\mathbb N}
\newcommand{\carn}{\hfill\rule{0.25cm}{0.25cm}}

\title[standard right factor of a Lyndon
  word]{Limit law of the length of the standard right factor of a
Lyndon
  word}

\begin{document}

\author{R{\'e}gine Marchand}
\author{Elahe Zohoorian Azad}
\address{Institut Elie Cartan Nancy (math{\'e}matiques)\\
Universit{\'e} Henri Poincar{\'e} Nancy 1\\
Campus Scientifique, BP 239 \\
54506 Vandoeuvre-l{\`e}s-Nancy  Cedex France}
\email{Regine.Marchand@iecn.u-nancy.fr, Elahe.Zohoorian@iecn.u-nancy.fr}

\subjclass[2000]{68R15,60B10,68Q25}
\keywords{random word, Lyndon word, standard right factor,
longest run, convergence in distribution}


\begin{abstract}
Consider the set of finite words on a totally ordered alphabet 
with $q$ letters. We prove that the distribution of the length of the standard right
factor of a random Lyndon word  with length $n$, divided by
$n$, converges to:
 $$\mu(dx)=\frac1q \delta_{1}(dx) + \frac{q-1}q \mathbf{1}_{[0,1)}(x)dx,$$
when $n$ goes to infinity. The convergence of all moments follows.
This paper completes thus the results of~\cite{Bassino}, giving the
asymptotics of the  mean length of the standard right factor of a random Lyndon word
with length $n$ in the case of a two letters alphabet.
\end{abstract}

\maketitle


\section{Introduction}
\label{intro}

Consider a finite totally ordered alphabet $\mathcal A$ and for each $n \in
\mathbb{N}=\{1,2,3,\dots\}$, denote by $\mathcal A^n$
the set of words with length $n$ on $\mathcal A$. A Lyndon word
with length $n$ is a word in $\mathcal A^n$ which is strictly
smaller, for the lexicographic order, than each of its proper
suffixes. We denote by $\mathcal L_n$ the set of Lyndon words with
length $n$.

\begin{figure}[t]
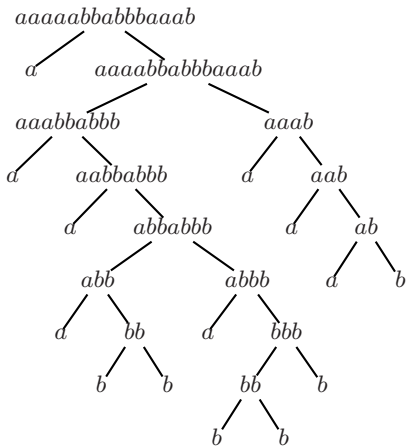

\label{olebodessin} {\small \arbre{$aaaaabbabbbaaab$}{
  \noeud{$a$}
  \arbre{$aaaabbabbbaaab$}{
      \arbre{$aaabbabbb$}{
          \noeud{$a$}
          \arbre{$aabbabbb$}{
              \noeud{$a$}
              \arbre{$abbabbb$}{
                  \arbre{$abb$}{
                      \noeud{$a$}
                      \arbre{$bb$}{
                          \noeud{$b$}
                          \noeud{$b$}
                      }
                  }
                  \arbre{$abbb$}{
                      \noeud{$a$}
                      \arbre{$bbb$}{
                          \arbre{$bb$}{
                              \noeud{$b$}
                              \noeud{$b$}
                          }
                          \noeud{$b$}
                      }
                  }
              }
          }
      }
      \arbre{$aaab$}{
          \noeud{$a$}
          \arbre{$aab$}{
              \noeud{$a$}
              \arbre{$ab$}{
                  \noeud{$a$}
                  \noeud{$b$}
              }
          }
      }
  }
} } \caption{Binary tree associated to the standard  decomposition
of the Lyndon word $aaaaabbabbbaaab$}
\end{figure}

The standard right factor $v$ of a Lyndon word $w$ is its smallest
proper suffix for the lexicographic order. Any Lyndon word $w$ can be written $uv$, in which
$u$ is a Lyndon word and $v$ is the standard right factor of $w$.
We call $uv$ the standard factorization of the Lyndon word $w$.
Lyndon words were introduced by Lyndon~\cite{Lyndon}, to
build a base of the free Lie algebra over $\mathcal A$. The
standard factorization plays a central part in the building
algorithm of this base. For each Lyndon word $w$, we can build a
binary tree in the following manner: say that $w$ is the root, and
has two children, that are the factors $u$ and $v$ of the standard
factorization of $w$. Since $u$ and $v$ are still Lyndon words,
they can also be divided into two standard factors which are their
children and so on (see figure \ref{olebodessin}). Then the
average height of these trees characterizes the complexity of the
building algorithm (see Chen, Fox and Lyndon~\cite{Chen} and Lothaire~\cite{Lothaire}). Thus the
informations on the length of the standard right factor of a
random Lyndon word are essential for the analysis of the
building algorithm.

For any Lyndon word $w \in \mathcal L_n$, let $R_{n}(w)$ denote the length of its
standard right factor. Endowing $\mathcal L_n$ with the
uniform probability measure
makes $R_n$
a random variable on $\mathcal L_n$. Bassino,
Cl{\'e}ment and Nicaud~\cite{Bassino}, with the help of generating functions,
 prove that the expectation $\E(R_n)$, in the case of a two letters alphabet, is asymptotically equal to $3n/4$. The aim of this
paper is to determine the limit distribution of $R_n/n$ as $n$
goes to infinity:

\begin{theorem}
\label{loilimiteq} For a totally ordered alphabet with $q$ letters, the
normalized length $R_n/n$ of a random Lyndon word of length $n$ converges in
distribution, when $n$ goes to infinity, to
$$\mu(dx)= \frac1q \delta_{1}(dx) + \frac{q-1}q \mathbf{1}_{[0,1)}(x)dx,$$
where $\delta_1$ denotes the Dirac mass on the point $1$, and $dx$
Lebesgue's measure on $\R$. All the moments of $R_n/n$ also
converge to the corresponding moments of the limit distribution.
\end{theorem}

\vspace{0.2cm} \noindent \textbf{Remark.}  In the case $q=2$, this
result was conjectured by Bassino, Cl{\'e}ment and
Nicaud~\cite{Bassino}. Simulations were also provided in this paper. 

\vspace{0.2cm} For the proof of this result, we focus first on the case of
a two letters alphabet, and then we indicate the way to adapt the proof for the case of
$q$ letters.

\vspace{0.2cm} Random Lyndon words are, in some sense,
conditioned random words. In section~\ref{words}, we obtain the number
of Lyndon words with length $n$ by dividing the number of primitive
words with length $n$ by $n$, (the shepherd's principle: counting the legs and
dividing by four to obtain the number of sheeps). Thus the
typical statistical behavior of a random word and of a random
Lyndon word can be easily linked (see lemma~\ref{usarg}).

\vspace{0.2cm} Our analysis starts in section~\ref{run}: we recall, among a
number of well known properties of random words with length $n$, those useful for our
purposes. In particular, we study the number of runs 
and the length of
the longest run of ``$a$'', which is typically of order $\log_2 n$. 

\vspace{0.2cm} The key step is to prove that the two longest runs
of ``$a$'' are approximately located along the word as two
independent uniform random variables, and thus the distance
$D_{n}$ between the first longest one and the second longest one
follows approximately
 the uniform law on $[0..n]$. The distance $D_{n}$ is of course closely
 related to the length of the
 standard
 right factor. We can distinguish two cases:

\begin{itemize}
\item
If the word obtained from the Lyndon word by deleting the first ``$a$'' is
still a Lyndon word, then the length of the standard right factor is equal
to $n-1 \simeq n$, and this happens with a probability close to $1/2$
(this probability is replaced by $1/q$ in the case of an alphabet with $q$ letters);
\item
Otherwise, the standard right factor begins by the  second longest run of
``$a$''. In this case,
the length of the standard right factor is equal to $n$ minus the
distance $D_{n}$, and is then approximately uniformly distributed
on $[0..n]$.
\end{itemize}

To prove that $D_{n}$ is approximately uniformly distributed
on $[0..n]$, we cut a random word with length
$n$ into distinct ``long blocks'' with length of order $\log_2 n$
(section~\ref{building}), in such a way that the long runs of ``$a$'' are at
the beginnings of the long blocks. Then we prove that the uniform distribution on Lyndon
words is invariant under  uniform permutation
of these blocks (section~\ref{permutation}). Thus the positions of the two smallest (for the lexicographic order) long blocks are approximately
uniformly distributed among all the possible positions of the
long blocks. As $n$ goes to infinity, the number of long blocks tends to infinity
and their lengths are negligible when compared to $n$. This leads to our
main  result, Theorem~\ref{distance}, which says
that the distance between the two smallest (for the lexicographic order) long blocks, divided
by the length $n$ of the word, follows asymptotically the uniform
law on $[0,1]$. In section~\ref{facteur}, we rephrase this result in terms
of standard right factor and finally we generalize, in section~\ref{section_q}, the obtained
results to the case of an alphabet with $q$ letters.


\section{Random words and random Lyndon words}
\label{words}


Let $\mathcal A=\{a,b\}$ be an ordered alphabet ($a<b$) and
 $\mathcal A^n$ be the set of all words with length $n$.
If $w \in \mathcal A^n$, write $w=(w_1, \dots, w_n)$ and define:
$$\tau w=(w_2,\dots, w_n,w_1).$$
Then $<\tau>=\{Id,\tau, \dots , \tau^{n-1}\}$ is the group of
cyclic permutations of the letters of a word with length $n$. A
word $w \in \mathcal A^n$ is called \textit{primitive} if 
$$(\exists k \in \{0,1, \dots, n-1\} \mbox{ such that }\tau^kw=w)
\Rightarrow (k
=0).$$
Denote by
$\mathcal P_n$ the set of primitive words in $\mathcal A^n$, by
$\mathcal N_n$ its complement. 

Remember that a Lyndon word
with length $n$ is a word in $\mathcal A^n$ which is strictly
smaller, for the lexicographic order, than each of its proper
suffix: it is equivalent to say that a word $w$ with length $n$ is a Lyndon word if and only if 
it is strictly smaller for the lexicographic order than every $\tau^kw$
with $k
\in \{1, \dots, n-1\}$.  We denote by $\mathcal L_n$ the set of Lyndon words with
length $n$.

The  group $<\tau>$ of cyclic permutations acts on $\mathcal A^n$,
and $\mathcal P_n$ and $\mathcal N_n$ are stable under this action.
Each orbit associated to a primitive word $w$ contains exactly
$n$ distinct words, and a unique Lyndon word, denoted by $\rho(w)$, which is the
smallest word in the orbit for the lexicographic order: the application $\rho$ is then the canonical projection of $\mathcal P_n$
on $\mathcal L_{n}$ associated to the action of $<\tau>$. 

\begin{exam} $\;$

$\bullet$ If $w=aabaaa$, then $\rho(w)=aaaaab$.

$\bullet$ If $A=\{aab,abb\}$ then
${{\rho}^{-1}}(A)=\{aab,aba,baa\}\cup\{abb,bab,bba\}$.
\end{exam}

\vspace{0.2cm} 
\noindent
As the set $\mathcal N_n$ of non-primitive words contains no Lyndon word, we
have, by the shepherd's principle, that:
$$card(\mathcal P_{n})= n \times card(\mathcal L_{n}).$$
Via the relation $card(\mathcal
P_n)=\sum_{d|n}{2^{n/d}\mu (d)},$ where $\mu$ is the M{\"o}bius function (see
the book by Lothaire~\cite{Lothaire}),
we are lead to:
$$card(\mathcal L_{n}) =
\frac{2^{n}}n(1+O(2^{-n/2}))\mbox{ and }card(\mathcal N_{n})=O(2^{n/2}).$$

 \vspace{0.2cm}
 In the sequel, we will consider  the two following probability spaces:
\begin{itemize}
\item
the set $\mathcal A^n$ of words with length $n$, endowed with the uniform probability~$\Pn$,
\item
the set $\mathcal L_n$ of Lyndon words with length $n$, endowed with the uniform probability $\Pnt$.
\end{itemize}
The probability measure $\Pnt$ can be seen as the conditional
probability on $\mathcal A^n$, given  $\mathcal L_n$.
The next lemma is obvious, but it is very useful in our proofs
because it allows to transfer results on random words to random
Lyndon words by neglecting non-primitive words and using the shepherd's
principle.

\begin{lemma}
\label{usarg} 
For $A \subset \mathcal L_n$, we have:
 $$\mid \tilde \Pn (A)- \P_{n}({\rho}^{-1} (A)) \mid \le O(2^{-n/2}).$$
\end{lemma}

 \noindent \textbf{Proof.} It is sufficient to note that 
$\Pnt(A)=\P_{n}({\rho}^{-1} (A)|\mathcal P_n)$ and 
$\Pn(\mathcal P_n)=1- O \left(2^{-n/2} \right).$


\section{Number of runs and length of the longest run}
\label{run}


This section deals with the number of runs and the length of the
longest run in a random Lyndon word. The
 results exposed in this section are not new, but are presented in
a convenient way for our proofs. The method is to get results for
random words, and to transfer them to random Lyndon words via lemma
\ref{usarg}.

\begin{defi}
Let $w$ be a word in $\mathcal A^n$. 
We denote by
$N_n(w)$ the number of runs in $w$, by $X_1(w),X_2(w), \dots,
X_{N_n}(w)$ their lengths, by $M_n(w)=\max \{ X_i(w), \; 1 \leq i
\leq  N_n(w)\}$ the length of the longest run in $w$ and by
$M_n^a(w)$ the length of the longest run of ``a'' in $w$.
\end{defi}

\begin{exam} $n=9$

\vspace{0.2cm}
\begin{tabular}{|c||c|c|c|c|}
\hline
& $N_6 $ & $(X_i)$ & $M_6$ & $M_6^a$ \\
\hline
\hline
$w=aabbbbaaa$ & $3$ & $2,4,3$ & $4$ & $3$ \\
\hline
$\rho(w)=aaaaabbbb$ & $2$ & $5,4$ & $5$ & $5$\\
\hline
\end{tabular}
\end{exam}

\vspace{0.1cm}
\begin{lemma}[Number of runs]
\label{Nn} For every $\gamma>0$, there
exists a constant $C_\gamma>0$ such that for any $\varepsilon>0$, 
\begin{eqnarray*}
\Pnt \left( \left| N_n-\frac{n}2 \right| \ge \gamma n^{-1/2+\varepsilon}
\right)
& \le &  O \left( \exp \left( -C_\gamma n^{2\varepsilon}
  \right) \right).
\end{eqnarray*}
\end{lemma}

\noindent \textbf{Proof.} First, we prove the above inequality for the
probability measure $\Pn$ on the set $\mathcal A^n$ of words with length
$n$. 
The cardinal of the event $\{N_{n}=k\}$ corresponds to the number of
compositions of 
the integer $n$ with $k$ parts (see Andrews~\cite{Andrews} and
Pitman~\cite{Pitman}):
$$\forall k \in \{1,2, \dots, n\}, \; \Pn (N_{n}=k)= \frac{1}{2^{n-1}}
\binom{n-1}{k-1}.$$ Thus $N_{n}-1$ is a binomial random variable
with parameters $(n-1,1/2)$, whose large deviations are well known
(see the book by Bollob{\'a}s~\cite[Th. 7, p.13]{Bollobas}  for instance): there exists a positive
constant $C'_\gamma$ such that
\begin{eqnarray*}
\Pn \left( \left| N_n-\frac{n}2 \right| \ge \gamma n^{-1/2+\varepsilon}
\right)
& \le &  O \left( \exp \left( -C_\gamma' n^{2\varepsilon}
  \right) \right).
\end{eqnarray*}

To obtain the same inequality for
the
probability measure $\Pnt$ on the set $\mathcal L_n$ of Lyndon words with
length $n$, note that for a primitive word $w$, we have $N_n(w) -1 \leq N_n(\rho(w))
\leq N_n(w)$. Thus, using Lemma \ref{usarg}, we obtain the announced
result. 
  \carn

\vspace{0.2cm} The next step is to study the length of the longest
run of a word $w \in \mathcal A^n$. For this mean, we will use the
following construction of the uniform probability measure on the set
of all infinite words on $\mathcal A$:

\begin{cons}
\label{cons2}
Let $(Z_i)_{i \in \N}$
be independent identically distributed geometrical random
variables with parameter $1/2$ defined on a probability space
$(\Omega, \mathcal F, \P)$ and let $\epsilon$ be a Bernoulli random
variable with parameter $1/2$ defined on $\Omega$ and independent of the $Z_i$'s.
To obtain a infinite random sequence of $a$ and
$b$, do the following:
\begin{itemize}
\item
if $\epsilon=1$, write $Z_1$ ``$a$'', followed by $Z_2$
``$b$'', followed by $Z_3$ ``$a$'' and so on...
\item
if $\epsilon=0$, write $Z_1$ ``$b$'', followed by $Z_2$
``$a$'', followed by $Z_3$ ``$b$'' and so on...
\end{itemize}
Truncating to keep the $n$ first letters gives a random variable
defined on $\Omega$ and uniformly distributed on $\mathcal A^n$.
Thus, in this setting, the number of runs is:
$$N_n(\omega)= \inf\left\{ k \in \N, \; \sum_{i=1}^k Z_i(\omega)\geq n
\right\},$$ and the lengths of runs are:
\begin{equation}
\label{EQruns} \forall i \in \{1, \dots, N_n-1\}, \; X_i=Z_i
\mbox{ and } X_{N_n}=n -\sum_{i=1}^{N_n-1}Z_i \leq Z_{N_n}.
\end{equation}
\end{cons}

We denote by $\log$ the natural logarithm, by $\log_2(x)=\log x/
\log 2$, by $\lceil x \rceil$ the smallest integer larger than $x$
and by $\lfloor x \rfloor$ the largest integer smaller than $x$.
The two next lemmas give estimates for the length of the longest
run in a random Lyndon word. These estimates
are related with the extreme values theory (see the books by Bingham,
Goldie and Teugels \cite{Bingham}, Resnick~\cite{Resnick} or the work of
Hitzenko and Louchard~\cite{Louchard}, or the initial works of Erd{\"o}s and
R{\'e}nyi~\cite{ErdosRenyi}
and Erd{\"o}s and R{\'e}v{\'e}sz~\cite{Erdos}).

\begin{lemma}[Longest run, small values]
\label{MnS}
For any  $\varepsilon>0$,
\begin{equation*}
\Pnt \left( M_n \le (1-\varepsilon)\log_2{n} \right) 
 \leq   O \left( \exp \left( -\frac{n^\varepsilon}4 \right)\right).
\end{equation*}
\end{lemma}

\noindent
 \textbf{Proof.} First, we prove the above inequality for the
probability measure $\Pn$ on the set $\mathcal A^n$ of words with length
$n$. The idea is that the number of the runs of
 ``$a$'' in a
 random word with length $n$ is highly concentrated around $n/4$ and that it is
 easy to estimate the maximum of $n/4$ independent geometrical random
 variables.

Note that $M_n \ge M_n^a$. Using (\ref{EQruns}) and lemma
\ref{Nn} with $\gamma=1$, we have:
\begin{eqnarray*}
&& \Pn(M_n \leq y ) \\
& \leq  & \Pn(M_n^a \leq y ) \\
& \leq  & \Pn \left( N_n \leq  \frac{n}{2} \left(
1-n^{-\frac12+\varepsilon} \right)
\right)
+ \Pn \left( X_i \leq y , \; 1 \le i \le  \frac{n}{4} \left(
1-n^{-\frac12+\varepsilon} \right) \right)    \\
& \leq & \Pn \left( N_n \leq  \frac{n}{2} \left(
1-n^{-1/2+\varepsilon} \right)
\right)
+ \Pn \left(Z_i \leq y, \;  1 \le i \le  \frac{n}{4} \left(
1-n^{-\frac12+\varepsilon} \right)  \right) \\
& \leq &  O \left( \exp \left( -C_1 n^{2\varepsilon}
  \right) \right) + (1-2^{\lfloor -y \rfloor})^{\lfloor \frac{n}{4} \left(
1-n^{-1/2+\varepsilon}  \right) \rfloor }
.
\end{eqnarray*}
To lighten notations, we consider the $\frac{n}{4} \left(
1-n^{-\frac12+\varepsilon} \right)$ first $Z_i$'s rather than the
$Z_i$'s corresponding to runs of ``$a$'', which would have obliged
us to distinguish whether the words begins with ``$a$'' or
``$b$''. Taking $y= (1 -\varepsilon)\log_2n$, we obtain
easily:
\begin{equation*}
\label{minmax} \Pn \left( M_n^a \leq (1 -\varepsilon) \log_2n
\right) \leq O \left( \exp \left( -\frac{n^\varepsilon}4 \right)
\right).
\end{equation*}

To obtain the same inequality for
the
probability measure $\Pnt$ on the set $\mathcal L_n$ of Lyndon words with
length $n$, note that, for a primitive word $w$, we have $M_n^a( w) \leq
 M_n^a(\rho(w))$. We can now use lemma \ref{usarg}
 to conclude. 
\carn

\begin{lemma}[Longest run, large values]
\label{lemMnGV}
For any  $1<A \le 2$,
\begin{equation*}
\Pnt \left( M_n^a \ge A\log_2{n} \right) \le \Pnt \left( M_n \ge A\log_2{n}
\right)
 \leq   O \left( n^{1-A} \right).
\end{equation*}
\end{lemma}

\noindent
 \textbf{Proof.} As before, we begin by proving the above inequality for the
probability measure $\Pn$ on the set $\mathcal A^n$. Note that we still have $M_n\ge M_n^a$.
 Thus for $y>0$, we have:
\begin{eqnarray}
\Pn(M_n^a \leq y )  
 & \geq & \Pn(M_n \leq y ) \nonumber \\
 & = & \Pn ( \forall i \in \{1,
\dots,N_n\}, \; X_i \leq
y) \nonumber \\
& \geq & \P (\forall i \in \{1, \dots,N_n\}, \; Z_i \leq y) \label{aaa}\\
& \geq & \P (\forall i \in \{1, \dots,n\}, \; Z_i \leq y) \label{bbb}\\
 & \geq &  (1-2^{- \lfloor y \rfloor})^n.\nonumber 
\end{eqnarray}
Inequality (\ref{aaa}) hold because of (\ref{EQruns}) and inequality
(\ref{bbb}) because $N_n \le n$. Taking $y= A\log_2n$, we obtain
easily the announced upper bound.

To come back to Lyndon words, note that for a primitive word $w$, we have
$M_n(\rho(w)) \leq \max\{M_n( w), X_1(w)+X_{N_n}(w)\}.$ Thus we
obtain:
\begin{eqnarray*}
&& \Pn \left(  M_n(\rho(w)) \geq  A \log_2{n} \right) \\
& \leq & \Pn \left( M_n(w)  \geq  A \log_2{n}\right) + \Pn \left(
X_1+X_{N_n} \ge A \log_2{n}  \right)\\
& \leq & \Pn \left( A \log_2{n}  \geq  M_n(w) \right) +\Pn \left(
 X_1 \ge  \frac{A}2   \log_2{n}
   \right)
  +\Pn \left( X_{N_n} \ge  \frac{A}2  \log_2{n}
   \right).
\end{eqnarray*}
Note that $X_1$ and $X_{N_n}$ have the same law. Thus using
(\ref{EQruns}), we get:  $$\Pn \left(
 X_1 \ge   \frac{A}2   \log_2{n}\right)
\le \Pn \left( Z_1 \ge   \frac{A}2 
  \log_2{n}\right)=O(n^{-\frac{A}2})=o\left(n^{1-A}\right),$$  by the
choice we made for $A$. Finally, by using Lemma
\ref{usarg}, we get the desired result. \carn


\section{Building of long blocks and short blocks of a word}
\label{building}


Let $0<\varepsilon<1$ and $B>2$ be fixed in this section. Our aim
here is to find, in a word $w \in \mathcal A^n$, some \emph{long blocks} beginning by a long run
of ``a'' in a word $w \in \mathcal A^n$, and we moreover want to choose them long enough to be
distinct with high probability. We study then the
positions of these long blocks along  the word $w$. Here is our definition: 
\begin{defi}[Long blocks]
Let $w$ be a word with length $n$. The long
blocks of $w$ are the subwords of $w$ that: 
\begin{itemize}
\item begin with a run of ``a'' with length equal or greater than
$(1-\varepsilon)\log_2 n$,
\item end with a run of ``b'' 
\item have the smallest possible length larger than $3\log_2n$.
\end{itemize}
We denote by $H_n$ the number of
long blocks.
\end{defi}

The next lemma estimates the number of long blocks for a random Lyndon
word. Note that although the crude estimate we give could be 
sharpen, it is sufficient for our mean.
\begin{lemma}[Number of long blocks]
\label{lemHn}
There exists a constant $D>0$ such that
\begin{eqnarray*}
\Pnt \left( \frac14  n^\varepsilon \leq H_n \leq \frac94
n^\varepsilon
\right) & \ge & 1-O \left( \exp \left(-D
    n^{\varepsilon} \right) \right).
\end{eqnarray*}
\end{lemma}

\noindent \textbf{Proof.} We begin once again by proving the inequality for the
probability measure $\Pn$ on the set $\mathcal A^n$. Set, for $i \geq 1$, $B_{i} =
\mathbf{1}_{\{ Z_{i}\geq(1-\varepsilon)\log_{2}{n}\}}.$ Then
$(B_i)_{i \geq 1}$ are independent identically distributed
Bernoulli random variables with parameter $p_{n, \varepsilon}$, 
which satisfies $n^{\varepsilon-1} \leq p_{n, \varepsilon} \leq
2n^{\varepsilon-1}.$ Note that
\begin{eqnarray*}
\sum_{1 \le 2i-1 \le N_n-1} B_{2i-1}(w) \leq H_n (w) \leq \sum_{i=1}^{n}
B_i(w) && \mbox{ if } w_1=a, \\
\sum_{1 \le 2i \le N_n-1} B_{2i}(w) \leq H_n (w) \leq \sum_{i=1}^{n}
B_i(w) && \mbox{ if } w_1=b.
\end{eqnarray*}
Therefore, by large deviation results for sums of independent
Bernoulli random variables (see for instance the book by Bollobas~\cite[Th. 7, p.13]{Bollobas}), there exists $D_{1}>0$ such that:
\begin{eqnarray*}
\Pn \left( H_n \geq \frac94 n^\varepsilon \right) & \leq & \P
\left(
\sum_{i=1}^{n} B_i \geq \frac94 n^\varepsilon \right) \\
& \leq & \P \left(\sum_{i=1}^{n} B_i - n\E B_1 \geq
\frac14 n^{\varepsilon} \right)
\leq O \left( \exp \left( -D_{1} n^{\varepsilon}
  \right) \right).
\end{eqnarray*}
In the same manner, by looking only to the $B_i$'s with odd indices
(when the word begins with ``$a$'') or only to the $B_i$'s with
even indices (when the word begins with ``$b$'') and using lemma
\ref{Nn} in which $\gamma=1$, we obtain the existence of $D_{2}>0$
such that:
\begin{eqnarray*}
&& \Pn \left( H_n \leq \frac14 n^\varepsilon \right) \\
& \leq & \P \left( N_n < \frac{n}{2} \left(1-n^{-1/2+\varepsilon}
\right)
  \right) + \P \left( \sum_{i=1}^{n/2\left(1-n^{-1/2+\varepsilon}
\right)} B_i\leq \frac14  n^\varepsilon
  \right)\\
& \leq & O \left( \exp \left( -C_1 n^{2\varepsilon}
  \right) \right) + O \left(  \exp \left(
    - D_{2} n^{\varepsilon} \right) \right).
\end{eqnarray*}
This proves the lemma for random words.

 For random Lyndon words, note that if $w$ is a primitive word, then
$H_n(w) \leq H_n(\rho(w)) \leq H_n(w)+1$ and use lemma \ref{usarg}.
\carn
\vspace{0.2cm}

The length of the long blocks has been chosen large enough to
ensure that two long blocks are distinct with high probability:

\begin{lemma}[Inequality of long blocks]
\label{En}
 Denote by $E_n$ the event that a word with length $n$ has at least two
 equal disjoint subwords with
 length at least $3\log_2n$.
 Then:
$$
\Pnt(E_n) \leq
O \left(n^{-1} \right).$$
\end{lemma}

\noindent \textbf{Proof.} We begin as usual with random words. By counting the number of possible subwords with
length $3\log_2n$ and their possible positions, we have:
$$\Pn \left(E_n \right)
\leq O \left\{ n^2 (2^{3{\log_2n}})(2^{-3{\log_2 n}})^2
\right\}\leq O \left(n^{-1}\right).$$ Lemma \ref{usarg} gives the
same estimate for Lyndon words. \carn

\vspace{0.2cm} We also want that the long blocks do not overlap
with high probability, or, in other words,  that the beginnings of 
long blocks are far away enough with high probability. This is
ensured by the next lemma:

\begin{lemma}[Minimal distance between beginnings of long blocks]
\label{Dn} Let $D_n$ be the event that  there exist at least two
long blocks which begin at a
  distance less than $8\log_2n$. Then:
$$\Pnt(D_n) \leq
O(n^{-(1-2\varepsilon)}\log_2 n ).$$
\end{lemma}

\noindent \textbf{Proof.} As usual, we start with the case of random words:
\begin{eqnarray*}
\Pn(D_n) & \leq & \Pn \left( H_n \geq \frac94 n^\varepsilon
\right)
+ \Pn \left(D_n\cap \left\{ H_n<
\frac 94 n^\varepsilon\right\} 
\right) 
\end{eqnarray*}
Let us denote by $F_n$ the last event. On $F_n$, at least one of the $H_n$ subwords with length $8\log_2
n$ starting just after a run of ``$a$'' with length
at least $(1-\varepsilon)\log_2 n$ must admit a subword of ``$a$''
with length $(1-\varepsilon)\log_2 n$ (which is the beginning of the next
long block). By an estimate analogous to the one used
in the previous lemma, 
$$\Pn \left(F_n \right)
\leq O \left\{ \frac 94 n^\varepsilon 8\log_2n
2^{-(1-\varepsilon){\log_2 n}}
\right\}=O( n^{-(1-2\varepsilon)} \log_2 n ).$$
By lemma \ref{lemHn}, the first term is negligible, and  Lemma \ref{usarg} concludes for Lyndon words.

\vspace{0.2cm}
Now we consider the set of ``good'' Lyndon words that satisfy all the previous
properties:
\begin{defi}[Good Lyndon words]
Denote by $\mathcal G_n$ the set of 
Lyndon words $w$ satisfying the following conditions:
\begin{itemize}
\item the maximal run of ``$a$'' satisfies $(1-\varepsilon)
\log_2{n}  \leq M_n^a \leq 2 \log_2{n}$ 
\item the  maximal run
satisfies $(1-\varepsilon) \log_2{n}  \leq M_n \leq 2 \log_2{n}$
\item the number of long blocks satisfies $\displaystyle \frac14
n^\varepsilon \leq H_n \leq \frac94 n^\varepsilon  $; 
\item the
beginnings of  long blocks are at a distance  at least $8
\log_2n$, in the sense $\mathcal G_n \subset D_n^c$; 
\item the word $w$ has no
equal long blocks, in the sense $\mathcal G_n \subset E_n^c$.
\end{itemize}
\end{defi}
 Note
that on $\mathcal G_n$, the length of a long block is less than
$3\log_2n+2\times2\log_2n= 7 \log_2{n}$, and that there is no overlapping
between two long blocks . The next lemma ensures that a
large proportion of Lyndon words are good Lyndon words:

\begin{lemma}
\label{Gn}
For every $n$ large enough:
$$\Pnt(\mathcal G_n) \geq
1- O \left( n^{-(1-2\varepsilon)}\log_2 n \right).$$
\end{lemma}

\noindent \textbf{Proof.} Everything has been proved in the previous lemmas \ref{MnS},
\ref{lemMnGV} (in which $A=2$), \ref{lemHn},
\ref{En} and \ref{Dn}. \carn

\vspace{0.2cm} Now, we note that a good Lyndon word $w\in
\mathcal G_n $ begins with a long block, ends with a run of
``$b$'', and all portions between long blocks  begin
with a run of ``$a$'' and end with a run of ``$b$''. We can thus give the
following definition of short blocks:

\begin{defi}[Short blocks]
For a good Lyndon word $w\in
\mathcal G_n $, we cut each section stretching between two long blocks into
short blocks, made of two consecutive runs of ``$a$'' and ``$b$''(in
this order).
\end{defi}

Note that short blocks have length equal or smaller to $4\log_2 n$. 


\section{Permutations of blocks for good Lyndon words}
\label{permutation}


In the previous section, we have cut any good Lyndon word $w$ into blocks
beginning with a run of ``$a$'' and ending with a run of ``$b$'':
the long ones and and the short ones. The long ones correspond to 
long runs of ``$a$'', and the first long block (at the beginning of $w$) is, by definition of a
Lyndon word, the smallest block for the lexicographic order.
We are going to see that we can keep this first long block of $w$ at the beginning of the word
 and permute the other blocks,
without changing the distribution on the set of good Lyndon words. 

In the
following, ``short'' and ``long'' refer to the type of blocks, while
``small'' and ``large'' refer to the lexicographic order on words. 

\begin{defi}[Permutation of blocks for good Lyndon words] Consider
  $w \in \mathcal G_n$.

1. We denote by $K_n(w)$ the total number of blocks, long and short, of $w$.

2. We denote by $(Y_i(w))_{0 \le i \le K_n(w)-1}$ the blocks of $w$ in
their order of appearance along $w$.
Certainly, the first block $Y_0(w)$ is the smallest block  among all blocks of~$w$.

3. Let $j_0(w)$ be the index of the second smallest block of $w$. 

4. We denote by $\mathfrak{S}_{K_n(w)-1}$ the set of permutations
of $\{1, \dots, K_n(w)-1\}$, and define
$$\sigma.w=Y_0(w)Y_{\sigma(1)}(w)\dots
Y_{\sigma(K_n(w)-1)}(w),$$ for $\sigma \in
\mathfrak{S}_{K_n(w)-1}$. Obviously, $\sigma.w \in \mathcal G_n$.

5. We define also $C(w)=\{\sigma.w, \sigma \in
\mathfrak{S}_{K_n(w)-1} \}$, the set of all words which are
obtained by the all the permutations of the blocks of $w$.
\end{defi}

The two cases of right factor exposed in the introduction can be rephrased
in the following manner: either  the standard
right factor is obtained by deleting the first ``$a$'', or it 
begins by the second smallest block of $w$, $Y_{j_0(w)}$.

In this section, we study the asymptotics of the position between the two
smallest blocks, and we will rephrase this result in terms of standard
right factor in the next section. 

Our main tool is the immediate following property: 

\begin{lemma}[Invariance in law under the permutations of blocks]  Let
$w_0$ be a fixed good Lyndon word. Consider the set
$\mathfrak{S}_{K_n(w_0)-1} \times
  C(w_0)$, endowed with the uniform probability. Then the random
variable
$$ \begin{array}{rrll}
W_{w_0}: & \mathfrak{S}_{K_n(w_0)-1} \times
  C(w_0) & \longrightarrow & C(w_0) \\
& (\sigma, w) & \longmapsto & \sigma.w
\end{array}
$$
follows the uniform law on $C(w_0)$.
\end{lemma}

\noindent \textbf{Proof.}
It is sufficient to note that, by construction, each word in  $C(w_0)$ has the same family of
blocks. 
\carn

\vspace{0.2cm} Thus, roughly speaking, the second smallest
block $Y_{j_0}$ has the same
probability to be at every possible place among all the blocks,
and this is why its position along the word $w$, divided by $n$, should follows approximately the
uniform
  law on $[0,1]$. To formalize this intuition and to exploit this
  invariance property, we enlarge our probability space $\mathcal G_n $:
  consider a sequence $(U_{i})_{i \in \N}$ of independent identically
  distributed random variables on a probability
space $(X, \mathcal X , \Q)$, following the uniform distribution on $[0,1]$. We denote by $\Pnb$ the uniform
probability on $\mathcal G_n $ and consider the product probability $\Q \otimes \Pnb $ on the
product space $X \times \mathcal G_n$; this means that
$(U_{i})_{i \in \N}$ are independent of the choice of the random
Lyndon word in $\mathcal G_n$.

\begin{defi}[Random permutation]
For $x \in X$, we define a uniform random permutation $\pi_x \in \mathfrak{S}_{K_n(w)-1}$ by the order statistics of
$U_{1}(x),U_{2}(x),...,U_{K_n(w)-1}(x)$:
$$U_{\pi_x(1)}(x) <
{U}_{\pi_x(2)}(x) < \dots < {U}_{\pi_x(K_n(w)-1)}(x).
$$
\end{defi}

Therefore from the previous lemma, the random variable
$$W: \;
\left\{ \begin{array}{rll}
X \times \mathcal G_n &
\longrightarrow & \mathcal G_n \\
(x, w) & \longmapsto & \pi_x.w
\end{array}
\right.
$$
follows the uniform law on $\mathcal G_n$.
We can now study, under the uniform probability $\Pnb$ on $\mathcal G_n$,
the position of the second smallest block defined by
$$d_n(w)=\frac1n \sum_{i=0}^{j_0(w)-1} |Y_i(w)|.$$ 
Here $|w|$
denotes the length of the word $w$.
Thus, the random variable
$$ \begin{array}{rrll}
d_n o W: & X \times  \mathcal G_n  &
\longrightarrow &
[0,1] \\
& (x, w) & \longmapsto & d_n(\pi_x.w)
\end{array}
$$
has the same law as the random variable $d_n$ under $\Pnb $.

We will thus focus on this new random variable to use the property
of invariance under the permutation of blocks. Remember
that the convergence in $L^2$ implies the convergence in
probability; thus, in the following, the notation $\|.\|_{\P}$ will denotes
the $L^2$-norm associated to a probability measure $\P$.

\begin{theorem}[Position of the second smallest
block] \label{distance} 
We~have:
$$ \|d_no W-U_{j_0}\|_{\Q \otimes \Pnb} \le
O \left( \left[ \frac{\log_2n}{n} \right]^{1/2}\right).
$$
This implies in particular that the law of
$d_n$ under $\Pnb$ converges weakly to the uniform law on $[0,1]$
and that every moment of $d_n$ converges to the corresponding moment
of the uniform distribution.
\end{theorem}

\noindent \textbf{Remark.} Coming back to random words,
this result implies that the normalized distance between the two smallest
blocks (which 
roughly corresponds to the two largest runs of ``$a$'') asymptotically
follows the uniform law on $[0,1]$.

\vspace{0.2cm}
\noindent \textbf{Proof.} We have:
\begin{eqnarray*}
d_n o W(x,w)=d_n(\pi_x.w) & = & \frac1n \left( |Y_0(w)|+\sum_{
j<\pi_x^{-1}(j_0(w))} |Y_{\pi_x(j)}(w)|
\right)\\
& = & \left( \frac1n |Y_0(w)| + \frac1n \sum_{j=1}^{K_n(w)-1}
|Y_j(w)| \mathbf{1}_{ \{U_{j}(x)<U_{j_0(w)}(x)\}} \right).
\end{eqnarray*}
By conditioning on $w$ and $U_{j_0(w)}$, and using the fact
that $\sum_i |Y_i|=n$, we
obtain:
\begin{eqnarray*}
&& \E_{\Q \otimes \Pnb} \left( \left. d_n(\pi_x.w) \right|
  w,U_{j_0(w)}\right)\\
& = & \left( \frac1n |Y_0(w)| + \frac1n 
\sum_{ \scriptsize
\begin{array}{cc}
1 \le j \le {K_n(w)-1}\\
j \neq j_0(w)
\end{array}}
|Y_j(w)| U_{j_0(w)} \right)\\
& = & \left( U_{j_0(w)} +(1-U_{j_0(w)})\frac{|Y_0(w)|}n -U_{j_0(w)}\frac{|Y_{j_0(w)}(w)|}n\right).
\end{eqnarray*}
On $\mathcal G_n$, $|Y_0|$ and $|Y_{j_0}|$ are bounded by $7 \log_2
n$, so 
\begin{eqnarray*}
&&\left\| \E_{\Q \otimes \Pnb} \left( \left.
      d_n(\pi_x.w) \right|
    w,U_{j_0(w)}\right) -U_{j_0(w)} \right\|_{\Q \otimes \Pnb}^2\\
&= &\E_{\Q\otimes \Pnb } \left[ \E_{\Q \otimes \Pnb} \left( \left.
      d_n(\pi_x.w) \right|
    w,U_{j_0(w)}\right) -U_{j_0(w)} \right]^2\\
&=& \E_{\Q \otimes \Pnb} \left[ \left( 1-U_{j_0(w)}
\right)\frac{|Y_0(w)|}n -U_{j_0(w)}\frac{|Y_{j_0(w)}(w)|}n\right]^2 \\
& \leq & \left[ \frac{14 \log_2n}{n} \right]^2 ,
\end{eqnarray*}
which tends to $0$ when $n$ goes to infinity. Now,
\begin{eqnarray}
&&\left\| d_n(\pi_x.w)  -\E_{\Q \otimes \Pnb} \left( \left.
      d_n(\pi_x.w) \right|
    w,U_{j_0(w)}\right) \right\|_{\Q \otimes \Pnb}^2\nonumber \\
&=&\E_{\Q \otimes \Pnb} \left[ d_n(\pi_x.w)  - \E_{\Q \otimes
\Pnb} \left( \left. d_n(\pi_x.w) \right| w,U_{j_0(w)}\right)
\right]^2 \nonumber \\
& = & \E_{\Q \otimes \Pnb} \left[  
\sum_{j=1}^{K_n(w)-1} \frac{|Y_j(w)|}n \mathbf{1}_{
\{U_{j}<U_{j_0(w)}\}} - \sum_{\scriptsize
\begin{array}{cc}
1 \le j \le {K_n(w)-1} \nonumber \\
j \neq j_0(w)
\end{array}}
\frac{|Y_j(w)|}n U_{j_0(w)} 
  \right]^2 \nonumber \\
& = & \E_{\Q \otimes \Pnb} \left[  
\sum_{\scriptsize
\begin{array}{cc}
1 \le j \le {K_n(w)-1} \nonumber \\
j \neq j_0(w)
\end{array}}
\frac{|Y_j(w)|}n \left( \mathbf{1}_{\{
U_{j}<U_{j_0(w)}\}} -U_{j_0(w)} \right)
 \right]^2 \nonumber \\
& \leq &\E_{\Q \otimes \Pnb} \left[  \frac14 \sum_{\scriptsize
\begin{array}{cc}
1 \le j \le {K_n(w)-1}\\
j \neq j_0(w)
\end{array}}
\left( \frac{|Y_j(w)|}{n} \right)^2
\right] \label{ccc}\\
& \leq & \frac{ 7\log_2n}{4n}, \label{ddd}
\end{eqnarray}
which tends to $0$ when $n$ goes to infinity. To obtain inequality
(\ref{ccc}), we conditioned first on $w$ and $U_{j_0(w)}$; for inequality
(\ref{ddd}), we used the facts that, 
on $\mathcal G_n$, all blocks have length smaller than
$7\log_2n$ and that  $\sum_i |Y_i|=n$.

Consequently,
\begin{eqnarray*}
\left\| d_n o W  - U_{j_0} \right\|_{\Q \otimes \Pnb} &
\leq & \left\| d_n(\pi_x.w)  -\E_{\Q \otimes \Pnb} \left( \left.
d_n(\pi_x.w) \right| w,U_{j_0}(w)\right)
\right\|_{\Q \otimes \Pnb}\\
&&+ \left\| \E_{\Q \otimes \Pnb} \left( \left. d_n(\pi_x.w)
\right|w,U_{j_0}(w)\right) -U_{j_0(w)}
\right\|_{\Q \otimes \Pnb} \\
& \leq & O \left( \left[ \frac{\log_2n}{n} \right]^{1/2}\right).
\end{eqnarray*}
The convergence of the other moments is a consequence of the
convergence in law, as $d_n$ is bounded by $1$.
\carn


\section{Limit distribution of the standard right factor}
\label{facteur}


In this section, we establish the convergence of the distribution of the normalized
length of the standard right factor of a random Lyndon word and give the
limit distribution, which follows quite easily from the result of the
previous section.
Remember that $\Pnt$ is the uniform probability on the set
$\mathcal L_n$ of Lyndon words with length $n$ and that $\Pnb$ is the uniform probability on the set
$\mathcal G_n$ of good Lyndon words with length $n$. The length of the
standard right factor of $w\in\mathcal L_{n}$ is denoted by
$R_n(w)$, and we introduce the normalized length of the standard
right factor $r_n(w)=R_n(w)/n$.

\begin{theorem}
\label{loilimite} 
As $n$ goes to infinity, $r_n$ converges in
distribution  to
$$\mu(dx)= \frac12 \delta_{1}(dx) + \frac12 \mathbf{1}_{[0,1)}(x)dx,$$
where $\delta_1$ denotes the Dirac mass at point $1$, and $dx$
Lebesgue's measure on $\R$. All the moments of $r_n$ also
converge to the corresponding moments of the limit distribution.
\end{theorem}

\vspace{0.2cm} 
\noindent 
\textbf{Proof.} First, we split the
set $\mathcal L_{n}$ in two parts, corresponding to the two cases of the introduction: 
\begin{itemize}
\item
$\mathcal L_{n}^1=a\mathcal
L_{n-1}$ contains exactly the Lyndon words $w$ whose standard
right factor is obtained by deleting the first ``$a$'' of the word
and has thus normalized length $r_n(w)=(n-1)/n$. Note
that
\begin{equation*}
\Pnt( \mathcal L_{n}^1)=\frac{card(\mathcal L_{n-1})}{card(\mathcal
  L_{n})}\sim \frac1{2}.
\end{equation*}
\item
$\mathcal
L_{n}^2= \mathcal L_{n} \backslash \mathcal L_{n}^1$ contains
exactly the Lyndon words $w$ whose standard right factor has
normalized length $r_n(w)$ strictly smaller than $(n-1)/n$. 
\end{itemize}

Now, forgetting the ``bad'' Lyndon words, using the inequality $r_n \le 1$
and lemma \ref{Gn}, we obtain the following inequality:
\begin{equation}
\label{EQU1}
\|r_n-r_n\mathbf{1}_{\mathcal G_n}\|_{\Pnt} \le (1-\Pnt(\mathcal G_n))^{1/2} \le O \left(
  n^{-\frac{1-2 \varepsilon}{2}} \sqrt{\log_2n} \right).
\end{equation}
But for $w \in \mathcal L_{n}^2 \cap \mathcal G_n$, the standard right
factor begins with the second smallest block $Y_{j_0}$ of $w$. Thus, in
this case, with the notations of the
previous section: $r_n(w)=1-d_n(w)$. Moreover, $\mathcal L_{n}^1
\cap \mathcal G_n$ and $\mathcal L_{n}^2 \cap \mathcal G_n$ are stable under the
permutations of blocks. Thus with the same setting as in
Theorem~\ref{distance},
\begin{eqnarray}
r_n(x,w) \mathbf{1}_{\mathcal G_n}(w)
&= &\left( (1-d_n(\pi_x.w)) \mathbf{1}_{\mathcal L_{n}^2 }(w) +
\frac{n-1}n \mathbf{1}_{\mathcal L_{n}^1}(w) \right) \mathbf{1}_{\mathcal G_n}(w)
\nonumber \\
&\stackrel{\mbox{\footnotesize{law}}}{=} & {r_n(w)}\mathbf{1}_{\mathcal G_n}(w) ,
\label{EQU3}
\end{eqnarray}
where the right hand side is a random variable from $ X \times
\mathcal L_n$, endowed with $\Q  \otimes \Pnt$ and the left hand
side is from $\mathcal L_n$, endowed with the uniform probability
$\Pnt$. Keeping in mind the result of the previous theorem, we
introduce for $(x,w) \in X \times \mathcal L_n$,
$$s_n(x,w) = \left( (1-U_{j_0(w)}(x)) \mathbf{1}_{\mathcal
    L_{n}^2 }(w) +\mathbf{1}_{\mathcal L_{n}^1}(w) \right)
\mathbf{1}_{\mathcal G_n}(w) .$$
Now, by Theorem \ref{distance},
\begin{eqnarray}
\label{EQU2}
&& \|r_n(x,w) \mathbf{1}_{\mathcal G_n}(w)-s_n(x,w)\|_{\Q  \otimes \Pnt} \nonumber \\
& \le & \|r_n(x,w) \mathbf{1}_{\mathcal G_n}(w)-s_n(x,w)\|_{\Q  \otimes \Pnb}
\nonumber \\
& \le & \left\| \left(r_n(x,w) -(1-U_{j_0(w)}(x)) \right)
  \mathbf{1}_{\mathcal
    L_{n}^2 }(w) -\frac1n\mathbf{1}_{\mathcal L_{n}^1}(w) \right\|_{\Q
\otimes
  \Pnb} \nonumber \\
& \leq & \left\| (d_n(x,w) - U_{j_0(w)}(x)) \right\|_{\Q  \otimes
\Pnb} + \frac1n \Pnb(\mathcal L_{n}^1)^{1/2} \nonumber \\
& \le & O \left( \left[ \frac{\log_2n}{n} \right]^{1/2}\right).
\end{eqnarray}
Note that the
position of $Y_{j_0}$, the second smallest block of $w$ in
$\pi_x.w$, is governed by $U_{j_0}$, which is clearly a uniform
random variable on $[0,1]$, independent of $w$. Note 
also that, thanks to lemma~\ref{Gn}, $\Pnb(\mathcal L_{n}^1 \cap \mathcal G_n) \sim 1/2$, and then
$\Pnb(\mathcal L_{n}^2 \cap \mathcal G_n) \sim 1/2$: consequently, the law of $s_n$
under $\Q  \otimes \Pnt$  converges weakly to $\mu(dx)= \frac12 [
\delta_{1}(dx) +
  \mathbf{1}_{[0,1]}(x)dx ]$. 

Now, as $\|r_n(x,w) \mathbf{1}_{\mathcal
    G_n}(w)-s_n(x,w)\|_{\Q  \otimes \Pnt}$ goes to $0$ by (\ref{EQU2}), a classical result
  (see for instance the book by Billingsley~\cite[Th. 4.2, p.25]{billingsley} in the first edition) ensures that the
  distribution of $r_n\mathbf{1}_{\mathcal
    G_n}$, as a random variable on $X\times \mathcal L_n$, also converges to $\mu(dx)= \frac12 [
\delta_{1}(dx) +
  \mathbf{1}_{[0,1]}(x)dx ]$. 
Using (\ref{EQU3}), the distribution of $r_n\mathbf{1}_{\mathcal
    G_n}$, as a random variable on $\mathcal L_n$, also converges to the
  same limit. Finally, (\ref{EQU1}) ensures the convergence of the
  distribution of $r_n$ to the
  same limit. 
\carn


\section{Generalization to the case of $q$ letters}
\label{section_q}

In this section, we generalize the previously obtained results to
the case of a totally ordered alphabet with $q$ letters: $\mathcal
A=\{a_{1},a_{2},\dots ,a_{q}\}$, $q\in\{2,3,4,...\}$
and $a_{1}< a_{2}<\dots<
a_{q}$. All the technics developed for
the simple case of two letters can be readily adapted in this
context and we just give the results and some indications for the
adaptations needed.

1. Denote by $\mathcal A^{n}$ the
set of words with length $n$ and by $\mathcal L_{n}=\mathcal L_{n}(\{a_{1},a_{2},\dots ,a_{q}\})$ the subset of Lyndon
words. The
probability measures $\Pn$ and $\Pnt$ are defined as before. As previously,
we have:
\begin{equation*}
\label{cardq}
 card(\mathcal L_{n})={\frac{q^{n}}{n}}(1+O(q^{-n/2})) .
\end{equation*}
The link between random Lyndon words and random words still holds:
 if $A \subset \mathcal L_n$, we have:
 $$\mid \tilde \Pn (A)- \P_{n}({\rho}^{-1} (A)) \mid \le O(q^{-n/2}).$$

2. Let $w$ be a word in $\mathcal A_{n}$. As previously, we
define its runs, its number of runs $N_{n}(w)$ and the length of
these runs $X_{1}(w), \dots ,X_{N_{n}}(w)$. To build these random
variables, we introduce a family $(Z_i)_{i \in \N}$ of
 independent identically distributed geometrical random
variables with parameter $(q-1)/q$ defined on a probability space
$(\Omega, \mathcal F, \P)$, and $(\epsilon_i)_{i \in \N}$, a family of
independent and identically distributed random variables with uniform distribution
on $\{a_{1},a_{2},...,a_{q}\}$, and independent of the $Z_i$'s. To obtain a random
sequence of letters, do the following:
\begin{itemize}
\item 
Select the letter $\epsilon_1$, and
write a run of $Z_1$ such letters. 
\item 
Select the letter $\epsilon_2$ conditioned to be distinct of
  $\epsilon_1$, and write a run of $Z_2$ such letters. 
\item 
Proceed by
recurrence: select the letter $\epsilon_{n+1}$
conditioned to be distinct of $\epsilon_n$, and
write a run of $Z_{n+1}$ such letters.
\end{itemize}
Truncating the $n$ first letters gives a random variable
defined on $\Omega$ and uniformly distributed on $\mathcal A_n$.
As previously,
$$N_n(\omega)= \inf\left\{ k \in \N, \; \sum_{i=1}^k Z_i(\omega)\geq n
\right\}.$$

3. Estimate the number of runs by using the fact that $N_n-1$
follows a binomial law with parameters $\left( n-1,\frac{q-1}q
\right)$ as in lemma~\ref{Nn}.

4. Estimate the length $M_n$ of the longest run and the length
$M_n^{a_1}$ of the largest run of $a_1$ as in lemma~\ref{MnS} and
lemma~\ref{lemMnGV} by using the same estimates on geometrical
laws. The typical order of $M_n$ and $M_n^{a_1}$ is $\log_q n$.

5. Define the long blocks:

\begin{defi}
Let $w$ be a word with length $n$. The long
blocks of $w$ are the subwords of $w$ that: 
\begin{itemize}
\item begin with a run of ``$a_1$'' with length equal or greater than
$(1-\varepsilon)\log_2 n$,
\item end just before an other  run of ``$a_1$'' (and consequently end with a
  run of a letter distinct from ``$a_1$''
\item have the smallest possible length larger than $3\log_q n$.
\end{itemize}
\end{defi}

Their number $H_n$ is, as in lemma~\ref{lemHn}, of order $n^\varepsilon$.
To prove this, introduce, for $i \geq 1$, the variable $B_{i} = \mathbf{1}_{\{
Z_{i}\geq(1-\varepsilon)\log_{q}{n}\}}.$
Then the $(B_i)_{i \geq 1}$ are independent identically distributed
Bernoulli random variables with parameter $p_{n, \varepsilon}$ satisfying
$n^{\varepsilon-1} \leq p_{n, \varepsilon} \leq
qn^{\varepsilon-1},$ thus we can have large deviation results.

6. We verify then that the long blocks do not overlap too often and
are distinct with high probability, as in lemmas~\ref{En} and
\ref{Dn}. Good Lyndon words are defined in the same manner as previously. Define the short blocks:
\begin{defi}
For a good Lyndon word $w\in
\mathcal G_n $, we cut each section stretching between two long blocks into
short blocks, that begin with a run of `$a_1$'' and end just before the
next run of ``$a_1$''.
\end{defi}

7. All is thus in place to permute the blocks as previously. With
the same setting as before, we obtain:

\begin{theorem}\label{distanceq} We have:
$$ \|d_no W-U_{j_0}\|_{\Q \otimes \Pnb} \le
O \left( \left[ \frac{\log_qn}{n} \right]^{1/2}\right).
$$
This implies in particular that the law of $d_n$ under $\Pnb$
converges weakly to the uniform law on $[0,1]$ and that every 
moment of $d_n$ converges to the corresponding moment of the limit law
\end{theorem}

8. To conclude for the length of the right factor, we split the
set $\mathcal L_{n}$ in two parts: 
\begin{itemize}
\item
$\mathcal
L_{n}^{1}=a_1\mathcal L_{n-1}(\{a_{1},a_{2},\dots ,a_{q}\}) \cup a_2\mathcal L_{n-1}(\{a_{2},\dots ,a_{q}\})
\cup \dots \cup a_{q-1}\mathcal L_{n-1}(\{a_{q-1},a_{q}\})$  contains exactly the
Lyndon words $w$ whose standard right factor is obtained by
deleting the first letter of the word and has thus normalized
length $r_n(w)=(n-1)/n$, 
\item
$\mathcal L_{n}^{q,2}= \mathcal
L_{n}^q \backslash \mathcal L_{n}^{q,1}$ contains exactly the
Lyndon words $w$ whose standard right factor has normalized length
$r_n(w)$ strictly smaller than $(n-1)/n$. 
\end{itemize}
The only difference is
that
\begin{equation*}
\Pnt( \mathcal L_{n}^{1})
 =  \frac{card(\mathcal L_{n-1}(\{a_{1},a_{2},\dots ,a_{q}\}))+
  \dots + card(\mathcal L_{n-1}(\{a_{q-1},a_{q}\}))}{card(\mathcal L_{n})}\\
 \sim  \frac1{q},
\end{equation*}
which gives Theorem \ref{loilimiteq}.


\section*{Acknowledgements}
We wish to thank Philippe Chassaing for pointing this problem to us,
and for many fruitful conversations.

\bibliographystyle{plain}
\bibliography{Lyndon}

\begin{thebibliography}{10}

\bibitem{Andrews}
G.~E. Andrews.
\newblock {\em The theory of partitions}.
\newblock Cambridge Mathematical Library. Cambridge University Press,
  Cambridge, 1998.
\newblock Reprint of the 1976 original.

\bibitem{Bassino}
F.~Bassino, J.~Clément, and C.~Nicaud.
\newblock The standard factorization of lyndon words: an average point of view.
\newblock {\em submitted to Elsevier Sciences, available at {\tt
  http://www-igm.univ-mlv.fr/\~\ bassino/biblio.html}}, 2003.

\bibitem{billingsley}
P.~Billingsley.
\newblock {\em Convergence of probability measures}.
\newblock Wiley Series in Probability and Statistics: Probability and
  Statistics. Second edition.

\bibitem{Bingham}
N.~H. Bingham, C.~M. Goldie, and J.~L. Teugels.
\newblock {\em Regular variation}, volume~27 of {\em Encyclopedia of
  Mathematics and its Applications}.
\newblock Cambridge University Press, Cambridge, 1989.

\bibitem{Bollobas}
B.~Bollob{\'a}s.
\newblock {\em Random graphs}, volume~73 of {\em Cambridge Studies in Advanced
  Mathematics}.
\newblock Cambridge University Press, Cambridge, second edition, 2001.

\bibitem{Chen}
K.-T. Chen, R.~H. Fox, and R.~C. Lyndon.
\newblock Free differential calculus. {IV}. {T}he quotient groups of the lower
  central series.
\newblock {\em Ann. of Math. (2)}, 68:81--95, 1958.

\bibitem{ErdosRenyi}
P.~Erd{\H{o}}s and A.~R{\'e}nyi.
\newblock On a new law of large numbers.
\newblock {\em J. Analyse Math.}, 23:103--111, 1970.

\bibitem{Erdos}
P.~Erd{\H{o}}s and P.~R{\'e}v{\'e}sz.
\newblock On the length of the longest head-run.
\newblock In {\em Topics in information theory (Second Colloq., Keszthely,
  1975)}, pages 219--228. Colloq. Math. Soc. J\'anos Bolyai, Vol. 16.
  North-Holland, Amsterdam, 1977.

\bibitem{Louchard}
P.~Hitzenko and G.~Louchard.
\newblock Distinctness of compositions of an integer: A probabilistic analysis.
\newblock {\em Random struct. Alg.}, 19:407--437, 2001.

\bibitem{Lothaire}
M.~Lothaire.
\newblock {\em Combinatorics on words}.
\newblock Cambridge Mathematical Library.

\bibitem{Lyndon}
R.~C. Lyndon.
\newblock On {B}urnside's problem.
\newblock {\em Trans. Amer. Math. Soc.}, 77:202--215, 1954.

\bibitem{Pitman}
J.~Pitman.
\newblock Combinatorial stochastic processes.
\newblock {\em Technical report}, 621, 2002.

\bibitem{Resnick}
S.~I. Resnick.
\newblock {\em Extreme values, regular variation, and point processes},
  volume~4 of {\em Applied Probability. A Series of the Applied Probability
  Trust}.
\newblock Springer-Verlag, New York, 1987.

\end{thebibliography}
\end{document}